\input amstex 
\input amsppt.sty 
\hsize 30pc
\vsize 43pc
\magnification=\magstep1
\def\nmb#1#2{#2}         
\def\idx{}               
\def\ign#1{}             

\redefine\o{\circ}
\define\X{\frak X}

\define\de{\delta}
\define\ep{\varepsilon}

\define\la{\lambda}
\define\rh{\rho}

\define\ph{\varphi}

\define\ps{\psi}

\redefine\i{^{-1}}
\define\row#1#2#3{#1_{#2},\ldots,#1_{#3}}
\define\x{\times}
\define\Fl{\operatorname{Fl}}
\define\ev{\operatorname{ev}}
\define\Mf{\Cal Mf}
\define\ddt{\left.\tfrac \partial{\partial t}\right\vert_0}
\redefine\L{{\Cal L}}
\def\today{\ifcase\month\or
 January\or February\or March\or April\or May\or June\or
 July\or August\or September\or October\or November\or December\fi
 \space\number\day, \number\year}
\topmatter
\title  Commutators of flows and fields
\endtitle
\author  Markus Mauhart \\
Peter W. Michor\footnote{Supported by Project P 7724 PHY 
of `Fonds zur F\"orderung der wissenschaftlichen Forschung'\hfill}
\endauthor
\leftheadtext{\smc Mauhart, Michor}
\affil
Institut f\"ur Mathematik, Universit\"at Wien, \\
Strudlhofgasse 4, A-1090 Wien, Austria.
\endaffil
\address
Institut f\"ur Mathematik, Universit\"at Wien,
Strudlhofgasse 4, A-1090 Wien, Austria
\endaddress
\email MICHOR\@AWIRAP.BITNET \endemail
\date \today   \enddate
\abstract The well known formula 
$[X,Y]=\tfrac12\tfrac{\partial^2}{\partial t^2}|_0
	(\Fl^Y_{-t}\o\Fl^X_{-t}\o\Fl^Y_t\o\Fl^X_t)$ for vector fields 
$X$, $Y$ is generalized to arbitrary bracket expressions and 
arbitrary curves of local diffeomorphisms.
\endabstract
\endtopmatter
\document

Let $M$ be a smooth manifold. It is well known that 
for vector fields $X,Y\in\X(M)$ we have
$$\align
0&=\tsize\ddt(\Fl^Y_{-t}\o\Fl^X_{-t}\o\Fl^Y_t\o\Fl^X_t),\\
[X,Y]&=\tfrac12\tfrac{\partial^2}{\partial t^2}|_0
	(\Fl^Y_{-t}\o\Fl^X_{-t}\o\Fl^Y_t\o\Fl^X_t).
\endalign$$
We give the following generalization:

\proclaim{\nmb.{1}. Theorem} Let $M$ be a manifold, let 
$\ph^i:\Bbb R\x M\supset U_{\ph^i}\to M$ be smooth mappings for 
$i=1,\dots,k$ where each $U_{\ph^i}$ is an open neighborhood of 
$\{0\}\x M$ in $\Bbb R\x M$, such that each $\ph^i_t$ is a 
diffeomorphism on its domain, $\ph^i_0=Id_M$, and 
$\ddt \ph^i_t=X_i\in\X(M)$. We put 
$[\ph^i_t,\ph^j_t]
	:=(\ph^j_t)\i\o(\ph^i_t)\i\o\ph^j_t\o\ph^i_t.$
Then for each formal bracket expression $B$ of length $k$ we have
$$\align
0&= \tfrac{\partial^\ell}{\partial t^\ell}|_0
	B(\ph^1_t,\dots,\ph^k_t)\quad\text{ for }1\le\ell<k,\\
B(X_1,\dots,X_k)&=\tfrac1{k!} \tfrac{\partial^k}{\partial t^k}|_0
	B(\ph^1_t,\dots,\ph^k_t)\in \X(M)
\endalign$$
in the sense explained in  \nmb!{3} below.
\endproclaim

In fact this theorem is a special case of the  more general theorem 
\nmb!{10} below. 
The somewhat unusual choice of the commutator of flows is explained 
by the fact that the bracket on the Lie algebra of the diffeomorphism 
group is the negative of the usual Lie bracket of vector fields.

\proclaim{\nmb.{2}. Lemma}
Let $c:\Bbb R\to M$ be a smooth curve. If 
$c(0)=x\in M$, $c'(0)=0,\dots,c^{(k-1)}(0)=0$, then $c^{(k)}(0)$ is a 
well defined tangent vector in $T_xM$ which is given by the 
derivation $f\mapsto (f\o c)^{(k)}(0)$ at $x$.
\endproclaim
 
\demo{Proof}
We have
$$\align
((f.g)\o c)^{(k)}(0)&= ((f\o c).(g\o c))^{(k)}(0) = \sum_{j=0}^k	
\tbinom kj (f\o c)^{(j)}(0)(g\o c)^{(k-j)}(0) \\
&= (f\o c)^{(k)}(0)g(x)+ f(x)(g\o c)^{(k)}(0),
\endalign$$
since all other summands vanish: $(f\o c)^{(j)}(0)=0$ for $1\le j 
<k$. \qed\enddemo

\subheading{\nmb.{3}. Curves of local diffeomorphisms}
Let $\ph :\Bbb R\x M\supset U_{\ph }\to M$ be a smooth mapping
where $U_{\ph }$ is an open neighborhood of 
$\{0\}\x M$ in $\Bbb R\x M$, such that each $\ph _t$ is a 
diffeomorphism on its domain and $\ph _0=Id_M$. We say that $\ph_t$ 
is a \idx{\it curve of local diffeomorphisms} though $Id_M$.

 From  lemma \nmb!{2} we see that if 
$\tfrac{\partial^j}{\partial t^j}|_0 \ph_t = 0$ 
for all $1\le j<k$, then 
$X:=\tfrac1{k!} \tfrac{\partial^k}{\partial t^k}|_0\ph_t$ 
is a well defined vector field on $M$. We say that $X$ is the 
\idx{\it first non-vanishing derivative} at 0 of the curve $\ph_t$ of 
local diffeomorphisms. We may paraphrase this as 
$(\partial^k_t|_0\ph_t^*)f=k!\L_Xf$. 

\subheading{\nmb.{4}. Natural vector bundles} See \cite{KMS, 6.14}.
Let $\Mf_m$ denote the category of all
smooth $m$-dimensional manifolds and local diffeomorphisms
between them. A \idx{\it vector bundle functor} or \idx{\it
natural vector bundle} is a functor $F$ which associates a
vector bundle $(F(M),p_M,M)$ to each manifold $M$ and a
vector bundle homomorphism 
$$\CD 
F(M) @>{F(f)}>> F(N) \\
@V{p_M}VV		 @VV{p_N}V \\
M    @>>{f}>    N     
\endCD$$
to each $f:M\to N$ in $\Mf_m$, which covers $f$ and is fiber wise a
linear isomorphism. If $f$ is the embedding of an open subset of $N$ 
then this diagram turns out to be a pullback diagram.
We also point out that $f\mapsto F(f)$ maps 
smoothly parameterized families to smoothly parameterized families, see
\cite{KMS, 14.8}. Assuming this property all vector bundle 
functors were classified by \cite{T}: They correspond to linear 
representations of  
higher jet groups, they are associated vector bundles to higher order 
frame bundles, see also \cite{KMS, 14.8}. 

Examples of vector bundle functors are
tangent and cotangent bundles, tensor bundles, and also the trivial 
bundle $M\x \Bbb R$ which will give us theorem \nmb!{1}.

\subheading{\nmb.{5}. Pullback of sections} Let $F$ be a vector bundle
functor on $\Mf_m$ as described in \nmb!{4}. Let $M$ be an 
$m$-manifold  and let $\ph_t$ be a curve of local diffeomorphisms 
through $Id_M$ on $M$. Then
the flow $\ph_t$, for fixed $t$, is a diffeomorphism defined
on an open subset $U_{\ph_t}$ of $M$.
The mapping
$$\CD
F(M) @<<< F(U_{\ph_t})  @>{F(\ph_t)}>>   F(M) \\
@V{p_M}VV  @VVV			@VV{p_M}V\\
M @<<< U_{\ph_t} @>>{\ph_t}>  M    
\endCD$$
is then a vector bundle isomorphism.

We consider a section $s\in C^\infty(F(M))$ of the vector
bundle $(F(M),p_M,M)$ and we define for $t\in \Bbb R$
$$\ph_t^*s := F(\ph_{t}\i)\o s\o \ph_t.$$
This is a local section of the bundle $F(M)$. For each $x\in M$
the value $(\ph_t^*s)(x)\in F(M)_x := p_M\i(x)$ is defined, if $t$ is
small enough. So in the vector space $F(M)_x$ the expression 
$\tfrac d{dt}|_0(\ph_t^*s)(x)$ makes sense and therefore the
section
$\tfrac d{dt}|_0(\ph_t)^*s$
is globally defined and is an element of $C^\infty(F(M))$. 
If $\ph_t=\Fl^X_t$ is the flow of a vector field $X$ on $M$  
this section 
$$\L_Xs := \tfrac d{dt}|_0(\Fl^X_t)^*s$$
is called the \idx{\it Lie derivative} of $s$ along $X$. It satisfies 
$\L_X\L_Y-\L_Y\L_X=\L_{[X,Y]}$, see \cite{KMS, 6.20}.

\proclaim{\nmb.{6}. Lemma} Let $\ph_t$ be a smooth curve of local 
diffeomorphisms through $Id_M$ with first non-vanishing derivative
${k!}X=\partial^k_t|_0\ph_t$. Then for any vector bundle 
functor $F$ and for any section $s\in C^\infty(F(M))$ we have the 
first non-vanishing derivative
$${k!}\L_Xs=\partial^k_t|_0\ph_t^*s.$$
\endproclaim

\demo{Proof}
This is again a local question, so let $x\in M$. We choose a complete 
Riemannian
metric on $M$ and we denote by $U_k$ the open ball with radius 
$r_k>0$ 
and center $x$ for this metric, and let let $\overline{U_k}$ be its 
closure. Since $\ph_0=Id_M$ we may choose a chart 
$(U,u:U\to \Bbb R^m)$ of $M$ with $x\in U$ and $u(U)=\Bbb R^m$, radii
$r_0>r_1>r_2>r_3>r_4>0$ and $\ep>0$ such that the following hold:
$\ph$ is defined and smooth on 
$((-2\ep,2\ep)\x U_0)$, $\ph([-\ep,\ep]\x \overline{U_1})\subset U$, 
$\ph((-\ep,\ep)\x U_2)\supset \overline{U_3}$,  
and $\ph((-\ep,\ep)\x \overline{U_4})\subset U_3$. 
Let $\Cal E$ be the set of all 
$f\in C^\infty(U_1,U)$ such that 
$f|\overline{U_2}$ is a diffeomorphism onto its image, $f(U_2)\supset 
\overline{U_3}$, and $f(\overline{U_4})\subset U_3$. 
Then via the linear isomorphism
$u_*:C^\infty(U_1,U)\to C^\infty(U_1,\Bbb R^m)$ which we suppress 
from now on,
the set $\Cal E$ is 
an open subset of the Frech\'et space 
$C^\infty(U_1,\Bbb R^m)$ for the  
compact $C^\infty$-topology, since the 
closures $\overline{U_k}$ are compact      
for each $r_k>0$ by completeness of the metric.

By cartesian closedness \cite{FK, 4.4.13} or 
\cite{KMb, 1.8} the curve 
$\check \ph:(-\ep,\ep)\to C^\infty(U_1,\Bbb R^m)$ is 
smooth and takes values in the open subset $\Cal E$.

\demo{Claim} Let $L(C^\infty(F(M)),C^\infty(F(U_4)))$ denote the 
space of all bounded linear mappings between the convenient vector 
spaces indicated which are equipped with the compact 
$C^\infty$-topology, and let
$P:C^\infty(U_1,\Bbb R^m)\supset \Cal E \to L(C^\infty(F(M)),C^\infty(F(U_4)))$ 
be the mapping given by $P(f)(s)=f^*s=F(f\i)\o s\o f$. Then $P$ is 
smooth. 
\enddemo

First we check that $P$ takes values in the space of bounded (i\. e\. 
smooth) linear mappings. We have to check that $P(f)$ maps smooth 
curves in $C^\infty(F(M))$ to smooth curves in $C^\infty(F(U_4))$.
A curve $c:\Bbb R\to C^\infty(F(M))$ is smooth if and only if the 
canonically associated mapping $\check c:\Bbb R\x M\to F(M)$ is 
smooth, see \cite{KMa,7.7.2}. 
But clearly $P(f)(c_t)(x)=(F((f|U_2)\i|U_3)\o c_t\o f|U_4)(x)$ 
is smooth in $(t,x)\in \Bbb R\x U_4$.

Now we check that $P$ itself is smooth, i\.e\. maps smooth curves in 
$\Cal E$ to smooth curves in $L(C^\infty(F(M)),C^\infty(F(U_4)))$.
So let $f:\Bbb R\to \Cal E\subset C^\infty(U_1,\Bbb R^m)$ be smooth, by 
cartesian closedness this means that 
$\hat f:\Bbb R\x U_1\to \Bbb R^m$ is 
smooth. By the finite dimensional implicit function theorem the 
mapping $(t,x)\mapsto f_t\i(x)$ is also smooth for 
$(t,x)\in \Bbb R\x U_3$. But then for each section 
$s\in C^\infty(F(M))$ the mapping
$(t,x)\mapsto (P(f_t)s)(x)= (F((f_t|U_2)\i|U_3)\o s\o f_t|U_4)(x)$
is also smooth since $F$ respects smoothly parameterized families.

By the smooth uniform boundedness 
principle \cite{FK, remark on page 89, also 4.4.7}, see also 
\cite{KMb, 1.7.2}, the assignment
$t\mapsto P(f_t)$ is smooth as a mapping
$$(-\ep,\ep)\to L(C^\infty(F(M)),C^\infty(F(U_4)))$$ 
if and only if the composition 
$$(-\ep,\ep)\to L(C^\infty(F(M)),C^\infty(F(U_4)))
	@>{\operatorname{ev}_s}>> 	C^\infty(F(U_4))$$
is smooth for each $s\in C^\infty(F(M))$. We have already checked 
this condition, so the claim follows. 

Now the smooth curve $\check\ph$ takes values in $\Cal E$, so we may 
compute for $1\le\ell\le k$ as follows:
$$\align
\partial^\ell_t|_0 \ph_t^*s 
	&= \partial^\ell_t|_0 (\ev_s\o P\o \check\ph)(t)\\
&= d(\ev_s\o P)(\ph_0)(\partial^\ell_t|_0\ph_t) + 0
\endalign$$
since each other term contains a derivative at 0 of $\ph_t$ of order 
less than $\ell$ which is 0, and thus we get
$\partial^\ell_t|_0 \ph_t^*s = 0$ for $\ell<k$ and 
$$\align
\partial^k_t|_0 \ph_t^*s &=
	d(\ev_s\o P)(Id_{U_2})(\partial_t^k|_0\ph_t)\\
&= 	d(\ev_s\o P)(Id_{U_2})(k!X)\\
&= k! d(\ev_s\o P)(Id_{U_2})(\partial_t|_0\Fl^X_t)\\   
&= k! \partial_t|_0 (\ev_s\o P\o \Fl^X)(t) = k! \partial_t|_0  
(\Fl^X_t)^*s = k!\L_Xs. \qed
\endalign$$
\enddemo

\proclaim{\nmb.{7}. Lemma} Let $M$ be a smooth manifold
and let $F$ be a vector bundle functor on $\Mf_m$.
Let $\ph_t$, $\ps_t$ be curves of local 
diffeomorphisms through $Id_M$ and let $s\in C^\infty(F(M))$ be a 
section of the vector bundle $F(M)\to M$.
Then we have 
$$\partial^k_t|_0(\ph_t\o\ps_t)^*s = 
	\partial^k_t|_0(\ps_t^*\ph_t^*)s = \sum_{j=0}^k\tbinom kj 
	(\partial^j_t|_0\ps_t^*)(\partial^{k-j}_t|_0\ph_t^*)s. $$
Also the multinomial version of this formula holds:
$$\partial^k_t|_0(\ph^1_t\o\dots\o\ph^\ell_t)^*s = 
	\sum_{j_1+\dots+j_\ell=k}\frac{k!}{j_1!\dots j_\ell!}
	(\partial^{j_\ell}_t|_0(\ph^\ell_t)^*)\dots
	(\partial^{j_1}_t|_0(\ph^1_t)^*)s.$$
\endproclaim

\demo{Proof} We only prove the binomial version.
The question is local on $M$, so let $U$ be an open neighborhood of 
some point $x$ in $M$ such that
$\ph$ is defined and smooth on 
$(-\ep,\ep)\x U$. 
 From the claim in the proof of lemma \nmb!{6} we know 
that $t\mapsto \ph_t^*$ is am smooth curve in the convenient vector space 
$L(C^\infty(F(M)),C^\infty(F(U)))$ of all bounded linear 
mappings.

Now let $V\subset M$ be an open neighborhood of $x$ 
such that $\ps$ is defined on 
$(-\ep,\ep)\x V$ and $\ps((-\ep,\ep)\x V)\subseteq U$. By the 
arguments just given the mapping $t\mapsto \ps_t^*$ is a smooth 
mapping $(-\ep,\ep)\to L(C^\infty(F(U)),C^\infty(F(V)))$ 
also. Composition 
$$\multline
L(C^\infty(F(M)),C^\infty(F(U)))\x 
	L(C^\infty(F(U)),C^\infty(F(V))) \to\\
	\to L(C^\infty(F(M)),C^\infty(F(V)))
\endmultline$$
is smooth and bilinear, see \cite{FK, 4.4.16} 
and we may just apply the Leibniz 
formula for higher derivatives of bilinear expressions of functions. 
We evaluate first at $s\in C^\infty(F(M))$ and then at $x\in M$ 
to obtain the formula
\qed\enddemo

\proclaim{\nmb.{8}. Lemma}
Let $\ph_t$ be a curve of local diffeomorphisms through 
$Id_M$ with first non-vanishing derivative 
${k!}X=\partial^k_t|_0\ph_t$. Then the inverse curve of local 
diffeomorphisms $\ph_t\i$ has first non-vanishing derivative 
$-{k!}X=\partial^k_t|_0\ph_t\i$.
\endproclaim

\demo{Proof}
For we have $\ph_t\i\o\ph_t=Id$, so by lemma \nmb!{7} we get for 
$1\le j\le k$
$$\multline
0=\partial^j_t|_0(\ph_t\i\o\ph_t)f = \sum_{i=0}^j\tbinom ji 
(\partial^i_t|_0\ph_t^*)(\partial^{j-i}_t(\ph_t\i)^*)f = \\
= \partial^j_t|_0\ph^*_t(\ph_0\i)^*f + 
\ph_0^*\partial^j_t|_0(\ph_t\i)^*f,
\endmultline$$
i\.e\. $\partial^j_t|_0\ph^*_tf= -\partial^j_t|_0(\ph_t\i)^*f$ as 
required.
\qed\enddemo

\proclaim{\nmb.{9}. Lemma} Let $M$ be a manifold, let  $F$ be a 
vector bundle functor, let $s$ be a smooth section of $F(M)$, 
let $\ph_t$ be a curve of local diffeomorphisms through 
$Id_M$ with first non-vanishing derivative 
${m!}X=\partial^m_t|_0\ph_t$, and let
$\ps_t$ be a curve of local diffeomorphisms through 
$Id_M$ with first non-vanishing derivative 
${n!}Y=\partial^n_t|_0\ps_t$.

Then the curve of local sections $[\ph_t,\ps_t]^*s = 
(\ps_t\i\o\ph_t\i\o\ps_t\o\ph_t)^*s$ has first non-vanishing derivative 
$${(m+n)!}\L_{[X,Y]}s = \partial^{m+n}_t|_0[\ph_t,\ps_t]^*s.$$
\endproclaim

\demo{Proof}
 From lemmas \nmb!{6} and \nmb!{8} we have the following first non-vanishing 
derivatives
$$\alignat2 {m!}\L_Xs &= \partial^m_t|_0\ph_t^*s,&\qquad
	{n!}\L_Ys &= \partial^n_t|_0\ps_t^*s,\tag1 \\
{m!}\L_{-X}s &= \partial^m_t|_0(\ph_t\i)^*s,&\qquad
	{n!}\L_{-Y}s &= \partial^n_t|_0(\ps_t\i)^*s.	
\endalignat$$
By the multinomial version of lemma \nmb!{7} we have
$$\align
A_Ns:&= \partial^N_t|_0(\ps_t\i\o\ph_t\i\o\ps_t\o\ph_t)^*s\\
&=\sum_{i+j+k+\ell=N}\frac{N!}{i!j!k!\ell!} (\partial^i_t|_0\ph_t^*)
	(\partial^j_t|_0\ps_t^*)(\partial^k_t|_0(\ph_t\i)^*) 
     (\partial^\ell_t|_0(\ps_t\i)^*)s.
\endalign$$
Let us suppose that $1\le n\le m$, the case $m\le n$ is similar.
If $N<n$ all summands are 0. If $N=n$ we have by lemma \nmb!{8}
$$A_Ns=(\partial^n_t|_0\ph_t^*)s + (\partial^n_t|_0\ps_t^*)s + 
	(\partial^n_t|_0(\ph_t\i)^*)s + (\partial^n_t|_0(\ps_t\i)^*)s = 
     0.$$
If $n<N\le m$ we have, using again lemma \nmb!{8}
$$\align
A_Ns&=\sum_{j+\ell=N}\frac{N!}{j!\ell!}
	(\partial^j_t|_0\ps_t^*)(\partial^\ell_t|_0(\ps_t\i)^*)s + 
     \de^m_N\left((\partial^m_t|_0\ph_t^*)s + 
     (\partial^m_t|_0(\ph_t\i)^*)s  \right)\\
&= (\partial^N_t|_0(\ps_t\i\o\ps_t)^*)s + 0 = 0.
\endalign$$
Now we come to the difficult case $m,n<N\le m+n$.
$$\align
A_Ns&= \partial^N_t|_0(\ps_t\i\o\ph_t\i\o\ps_t)^*s 
	+ \tbinom Nm (\partial^m_t|_0\ph_t^*)
	(\partial^{N-m}_t|_0(\ps_t\i\o\ph_t\i\o\ps_t)^*)s\\
&\quad + (\partial^N_t|_0\ph_t^*)s, \tag 2
\endalign$$
by lemma \nmb!{7}, since all other terms vanish, see \thetag4 below. By 
lemma \nmb!{7} again we get:
$$\align
\partial^N_t|_0&(\ps_t\i\o\ph_t\i\o\ps_t)^*s = 
	\sum_{j+k+\ell=N}\frac{N!}{j!k!\ell!}
	(\partial^j_t|_0\ps_t^*)(\partial^k_t|_0(\ph_t\i)^*) 
     (\partial^\ell_t|_0(\ps_t\i)^*)s \\
&=\sum_{j+\ell=N}\tbinom Nj (\partial^j_t|_0\ps_t^*)
	(\partial^\ell_t|_0(\ps_t\i)^*)s
	+ \tbinom Nm (\partial^{N-m}_t|_0\ps_t^*)
	(\partial^m_t|_0(\ph_t\i)^*)s \tag3\\
&\quad + \tbinom Nm (\partial^m_t|_0(\ph_t\i)^*)
	(\partial^{N-m}_t|_0(\ps_t\i)^*)s
	+ \partial^N_t|_0(\ph_t\i)^*s\\
&= 0+\tbinom Nm (\partial^{N-m}_t|_0\ps_t^*)
	m!\L_{-X}s + \tbinom Nm 
     m!\L_{-X}(\partial^{N-m}_t|_0(\ps_t\i)^*)s\\
&\quad + \partial^N_t|_0(\ph_t\i)^*s,\qquad\text{ using \thetag 1}\\
&= \de^N_{m+n}(m+n)!(\L_X\L_Y-\L_Y\L_X)s
	+ \partial^N_t|_0(\ph_t\i)^*s\\
&= \de^N_{m+n}(m+n)!\L_{[X,Y]}s
	+ \partial^N_t|_0(\ph_t\i)^*s
\endalign$$
 From the second expression in \thetag3 one can also read off that 
$$\partial^{N-m}_t|_0(\ps_t\i\o\ph_t\i\o\ps_t)^*s = 
	\partial^{N-m}_t|_0(\ph_t\i)^*s.\tag 4$$  
If we put \thetag3 and \thetag4 into \thetag2 we get, using lemmas 
\nmb!{7} and \nmb!{8}
again, the final result which proves lemma \nmb!{9}:
$$\align
A_Ns&= \de^N_{m+n}(m+n)!\L_{[X,Y]}s + \partial^N_t|_0(\ph_t\i)^*s\\
&\quad+\tbinom Nm(\partial^m_t|_0\ph_t^*)(\partial^{N-m}_t|_0(\ph_t\i)^*)s
	+ (\partial^N_t|_0\ph_t^*)s\\
&=\de^N_{m+n}(m+n)!\L_{[X,Y]}s + \partial^N_t|_0(\ph_t\i\o\ph_t)^*s\\
&=\de^N_{m+n}(m+n)!\L_{[X,Y]}s + 0.\qed
\endalign$$ 
\enddemo

\proclaim{\nmb.{10}. Theorem}
Let $M$ be a manifold, let 
$\ph^i$ be smooth curves of local diffeomorphisms through $Id_M$ for 
$i=1,\dots,j$ with non-vanishing first derivative
$\partial^{k_i}_t|_0 \ph^i_t=k_i!\,X_i\in\X(M)$.
Let $F$ be a vector bundle functor and let $s\in C^\infty(F(M))$	be a 
section.
Then for each formal bracket expression $B$ of length $j$ we have
$$\align
0&= \tfrac{\partial^\ell}{\partial t^\ell}|_0
	B(\ph^1_t,\dots\ph^k_t)^*s\quad\text{ for }1\le\ell<k,\\
\L_{B(X_1,\dots,X_k)}s&=\tfrac1{k!} \tfrac{\partial^k}{\partial t^k}|_0
	B(\ph^1_t,\dots\ph^k_t)^*s\in C^\infty(F(M)),
\endalign$$
where $k=k_1+\dots+k_j$.
\endproclaim

\demo{Proof}
Apply lemma \nmb!{9} recursively. 
\qed\enddemo

\proclaim{\nmb.{11}. Proposition} Let $\ph$ be a curve of local 
diffeomorphisms through $Id_M$ with first non-vanishing derivative
$k!X=\partial^k_t|_0\ph_t$. Then the curve of local vector
fields  $(\partial_t\ph_t)\o\ph_t\i$ has as first non-vanishing 
derivative
$$k!X=\partial^{k-1}_t|_0\left((\partial_t\ph_t)\o\ph_t\i\right).$$
\endproclaim

\demo{Proof}
Using lemma \nmb!{7} for $f\in C^\infty(M,\Bbb R)$ we have for 
$1\le\ell<k$:
$$\align
\partial^{\ell-1}_t|_0((\partial_t\ph_t)\o\ph_t\i)f &= 
\partial^{\ell-1}_t|_0(\ph_t\i)^*\partial_t\ph_t^*f\\
&= \sum_{j=0}^{\ell-1}\tbinom {\ell-1}j 
	(\partial^j_t|_0(\ph_t\i)^*)(\partial^{\ell-j}_t|_0\ph_t^*)f\\
&= (\ph_0\i)^*(\partial^{\ell}_t|_0\ph_t^*)f +0 \\
&= \de^k_\ell k!\L_Xf.\qed
\endalign$$
\enddemo

\proclaim{\nmb.{12}. Corollary} Let $G$ be a Lie group with Lie 
algebra $\frak g$. For $g$, $h\in G$ we consider the group commutator
$[g,h]=ghg\i h\i$. Then for any bracket expression $B$ of length $k$ 
and $X_i\in \frak g$ we have 
$$\multline
k!B(X_1,\dots,X_k) = \partial^k_t|_0B(\exp tX_1,\dots,\exp tX_k)  \\
= \partial^{k-1}_t|_0\left(T\la_{B(\exp tX_1,\dots,\exp tX_k)\i)}
	(\partial_tB(\exp tX_1,\dots,\exp tX_k))\right),
\endmultline$$
where $\la_g$ denotes left translation by $g$.
\endproclaim
The first equation is a generalization of the well known `Trotter product 
formula', i\. e\. the case of $B=[\quad,\quad]$.

\demo{Proof}
The flow of the left invariant vector field $L_X$ corresponding to 
$X\in\frak g$   
is the right translation $\rh_{\exp tX}$ by $\exp tX$, so we  
just apply theorem \nmb!{1} to get
$$\align
k!B(L_{X_1},\dots,L_{X_k}) &= 
	\partial^k_t|_0B(\rh_{\exp tX_1},\dots,\rh_{\exp tX_k})  \\
&= \partial^k_t|_0\rh(B(\exp tX_1,\dots,\exp tX_k)),		 \tag1
\endalign$$
where in the first line the commutator of flows is applied, and in 
the second line the group commutator with reversed order. Evaluating 
both sides at $e\in G$ gives the first formula.
 From \thetag1 and proposition \nmb!{11} we get
$$\multline
k!B(L_{X_1},\dots,L_{X_k}) = \\ 
= \partial^{k-1}_t|_0\left(
  	(\partial_t\rh(B(\exp tX_1,\dots,\exp tX_k)))\o
	\rh(B(\exp tX_1,\dots,\exp tX_k))\i\right).
\endmultline$$
We evaluate this at $e\in G$ and get
$$\align
k!B&(X_1,\dots,X_k) = \\ 
&= \partial^{k-1}_t|_0\left(
  	(\partial_t\rh(B(\exp tX_1,\dots,\exp tX_k)))
	(B(\exp tX_1,\dots,\exp tX_k)\i)\right)\\
&= \partial^{k-1}_t|_0\left(T\la_{B(\exp tX_1,\dots,\exp tX_k)\i)}
	(\partial_tB(\exp tX_1,\dots,\exp tX_k))\right).\qed
\endalign$$
\enddemo

\enddocument